\documentclass[twoside]{article}
\usepackage{graphicx,amssymb,mathrsfs,amsmath,color,fancyhdr,calc,hyperref,euscript}
\renewcommand{\baselinestretch}{1.06} 
\catcode`@=11 \long\def\@makefntext#1{\noindent #1}
\newskip\tabcentering \tabcentering=1000pt plus 1000pt minus 1000pt
\def\REF#1{\par\hangindent\parindent\indent\llap{#1\enspace}}%\ignorespaces}
\def\MCH#1#2{\setbox0=\hbox{\raise#1\hbox{#2}}\smash{\box0}}% move char

\def\dl{\displaystyle}
\let\@oddfoot\@empty  \let\@evenfoot\@empty

\def\@evenhead{}\def\@oddhead{}
\def\@evenhead{\vbox{\hbox to \textwidth{\footnotesize\rm\hbox to
1.0cm{\thepage\hfill} \hfill\hspace{2mm}\footnotesize{ \emph{Liu B
$\&$ Yu J}}}}}
\def\@oddhead{\vbox{\hbox to \textwidth{\footnotesize
{\it On the Anomaly Formula for the Cappell-Miller Holomorphic
Torsion} \hfill{\ } \hfill\hbox to 1cm{\hfill\thepage}}}}

\numberwithin{equation}{section} %公式按照section编号
\def\th#1{\vspace{1mm}\noindent{\bf #1}\quad }    %规定定理的格式
\usepackage{amsthm}
\newtheoremstyle{mystyle}{1mm}{0pt}{\it}{}{\bf}{.}{10pt}{}
\theoremstyle{mystyle}
\newtheorem{thm}{Theorem}[section]%定理按照section编号
\newtheorem{defn}[thm]{Definition}%定义
\newtheorem{prop}[thm]{Proposition}%命题
\newtheorem{lem}[thm]{Lemma}%引理

\def\pf#1{\vspace{1mm}\noindent{\it #1.}\quad}%输出“Proof.”

\def\leq{\leqslant}
\def\geq{\geqslant}
\def\no{\noindent}
\def\R{{\Bbb R}}       
\def\C{{\Bbb C}}  \def\hml{\end{document}}  \newsymbol\wjzhml 203F

\begin{document}
\abovedisplayskip=3pt plus 1pt minus 1pt %公式以上距离
\belowdisplayskip=3pt plus 1pt minus 1pt %公式以下距离

\def\le{\leqslant}
\def\ge{\geqslant}
\def\dl{\displaystyle}
\def\ch{{\rm  ch}}
\def\td{{\rm  Td}}
\def\Hom{{\rm  Hom}}
\def\End{{\rm  End}}
\def\tr{{\rm Tr}}
\newcommand{\Ker}{{\rm Ker}}
\newcommand{\wi}{\widetilde}

%Vol. 50 No. 1 75--85

%-------------------  First Head  -----------------------------------------

%===================Text=============================================

\vspace{8true mm}

\renewcommand{\baselinestretch}{1.9}\baselineskip 19pt
%\leftline{\Large\sf MATHEMATICS}

\noindent{\LARGE\bf On the Anomaly Formula for the Cappell-Miller Holomorphic Torsion}%题目

\vspace{0.5 true cm}

\noindent{\normalsize\sf Bo LIU$^{1\,\dag}$ \&  Jianqing YU$^{2}$%作者
\footnotetext{\baselineskip 10pt $^\dag$ Corresponding author%\\Thiswork was
} }

\vspace{0.2 true cm}
\renewcommand{\baselinestretch}{1.5}\baselineskip 12pt
\noindent{\footnotesize\rm $^1$ Chern Institute of Mathematics \&
LPMC, Nankai University, Tianjin 300071, P. R. China
\\
$^2$ Chern Institute of Mathematics \& LPMC, Nankai University,
Tianjin 300071, P. R. China
\\
(email: boliumath@mail.nankai.edu.cn, yjqxc@mail.nankai.edu.cn)\vspace{4mm}
}%联系方式

\baselineskip 12pt \renewcommand{\baselinestretch}{1.18}
\noindent{{\bf Abstract}\small\hspace{2.8mm} %摘要
We present an explicit expression of  the anomaly formula for the
Cappell-Miller holomorphic torsion for K\"{a}hler manifolds.}

\vspace{1mm} \no{\footnotesize{\bf Keywords:\hspace{2mm}
%关键词
Cappell-Miller holomorphic torsion, anomaly formula, heat kernel}}

\no{\footnotesize{\bf MSC(2000):\hspace{2mm}
58J52%分类号
}}

\vspace{2mm}
\baselineskip 15pt
\renewcommand{\baselinestretch}{1.22}
\parindent=10.8pt
\rm\normalsize\rm
%---------------------------------------------1------------------------------------------------------------------
\makeatletter
\@startsection{section}{0}{0pt}{2mm}{0.5mm}{\bfseries}
{Introduction}
\makeatother

As a subsequent paper of the celebrated work \cite{ray1971rta}, Ray
and Singer defined a holomorphic torsion associated to the
$\overline{\partial}$-complex of complex manifolds in
\cite{ray1973atc}. In \cite[(2.5)]{ray1973atc}, in the flat holomorphic case, they obtained the variation
of this holomorphic torsion with respect to the change of Hermitian metrics
of the underlying complex manifold, as the constant coefficient of the small
time asymptotic expansion of certain trace of the associated heat kernel. Further in
\cite{bismut1988ata3}, using probability method, Bismut, Gillet and
Soul\'{e} obtained the explicit anomaly formula for the case where the
holomorphic bundle is endowed with Hermitian metrics and the base manifold is assumed to be K\"ahler.

In the recent paper of Cappell and Miller \cite{cappell2010cvant},
the holomorphic torsion is extended to coupling with an arbitrary
holomorphic bundle with compatible connection of type (1,1).
However, comparing with the operators dealt with in
\cite{bismut1988ata3}, in the present general setting the associated
operators are not necessarily self-adjoint and the torsion is
complex valued. Especially, the torsion is independent of the given
Hermitian metric of the holomorphic bundle.

In this paper, we get an explicit expression of the anomaly formula for this
Cappell-Miller holomorphic torsion for K\"ahler manifolds. As is
obtained in \cite{cappell2010cvant}, the variation of the
holomorphic torsion is the constant term in the Laurent expansion of
the supertrace of a certain smooth kernel. To compute the constant
term, following \cite[Section 1(h)]{bismut1988ata3} in spirit, we
introduce the Grassmann variables ${\rm d}a,{\rm d}\overline{a}$ and
identify the constant term as the coefficient of ${\rm d}a{\rm
d}\overline{a}$ of $\lim_{t\rightarrow0^+}\tr_s[\exp(-tI_t')]$,
where for any $t>0$, $I_t'$ is a generalized non-self-adjoint
Laplacian with parameters ${\rm d}a,{\rm d}\overline{a}$.

In our final calculation, we don't use probability theory as in
\cite{bismut1988ata3}, but modify the proof which is presented in
\cite[Chapter 4]{berline2003hka} to deduce the local index theorem
\cite[Theorem 4.1]{berline2003hka}. The technical difficulties in the
modification are the singular terms turning up in the rescaled
operator and the convergence of the rescaled heat kernel for the
introduced rescaling parameter $\varepsilon$ as
$\varepsilon\rightarrow 0$. We first modify the estimates on heat
kernels given in \cite[Chapter 2]{berline2003hka}, which essentially
depend on the properties of the introduced Grassmann variables, then
using Donnelly's conjugation technique and the general  Mehler
formula for a generalized harmonic oscillator (cf.
\cite{donnelly1988lit}), we overcome the difficulties.

The rest of this paper is organized as follows. In section 2, we
recall the basic definition of the Cappell-Miller torsion of the
holomorphic bundle endowed with a compatible (1,1) connection over
complex manifolds and state the anomaly formula for the case where
the base manifold assumed to be K\"{a}hler. In section 3, we deduce
some estimates on the heat kernels in a little more general case
than what we need. In section 4, using the local index theorem
techniques, we prove the anomaly formula stated in section 2.

%------------------------------------------------2---------------------------------------------------------------

\makeatletter \@startsection{section}{0}{0pt}{2mm}{0.5mm}{\bfseries}
{Cappell-Miller torsion of the holomorphic bundle with a compatible
(1,1) connection}\makeatother

In this section, we recall the definition of the Cappell-Miller
torsion of the holomorphic bundle with a compatible (1,1) connection
and state an anomaly formula for the Cappell-Miller torsion under
the K\"{a}hler condition.

%-------------------------------2.1-------------
\makeatletter
\@startsection{subsection}{0}{0pt}{2mm}{0.2mm}{\bfseries} {The
Cappell-Miller holomorphic torsion} \makeatother

Let $(M,J)$ be a complex manifold with complex structure $J$ and its
complex dimension be $n$. Let $TM$ be the corresponding real tangent
bundle. Let $g^{TM}$ be any Riemannian metric on $TM$ compatible
with $J$ and $\nabla^{TM}$ be the corresponding Levi-Civita
connection.

Let $E\rightarrow M$ be a complex holomorphic bundle over $M$
endowed with a connection $\nabla^E$. Let $g^E$ be a Hermitian
metric on $E$.

For $0\leq r\leq 2n$, let
$\Omega^{r}(M,E)=\Gamma\left(M,\Lambda^r(T^*M)\otimes E\right)$ be
the space of smooth $r$-forms on $M$ with values in $E$.

The complex structure $J$ induces a splitting
$TM\otimes_{\R}\C=T^{(1,0)}M\oplus T^{(0,1)}M$, where $T^{(1,0)}M$
and $T^{(0,1)}M$ are the eigenbundles of $J$ corresponding to the
eigenvalues $\sqrt{-1}$ and $-\sqrt{-1}$, respectively. Let
$T^{*(1,0)}M$ and $T^{*(0,1)}M$ be the corresponding dual bundles.

For any
$0\leq p,q\leq n$, let
\begin{equation*}
\Omega^{p,q}(M,E)=\Gamma\left(M,\Lambda^p\bigl(T^{*(1,0)}M\bigr)\otimes\Lambda^q\bigl(T^{*(0,1)}M\bigr)\otimes
E\right)
\end{equation*}
be the space of smooth $(p,q)$-forms on $M$ with values in $E$. Set
\begin{equation*}
\Omega^{\ast,\ast}(M,E)=\bigoplus_{p,q=0}^{n}\Omega^{p,q}(M,E).
\end{equation*}

We have the direct sum decomposition $\Omega^r(M,\C)={\textstyle
\bigoplus}_{p+q=r}\Omega^{p,q}(M,\C)$ and the differentials ${\rm
d}:\Omega^r(M,\C)\rightarrow\Omega^{r+1}(M,\C)$,
$\partial:\Omega^{p,q}(M,\C)\rightarrow\Omega^{p+1,q}(M,\C)$ and
$\overline{\partial}:\Omega^{p,q}(M,\C)\rightarrow\Omega^{p,q+1}(M,\C)$
with ${\rm d}=\partial+\overline{\partial}$.

We will denote by $\langle\cdot,\cdot\rangle$ the $\C$-bilinear form
on $TM\otimes_{\R}\C$ induced by $g^{TM}$. $T^{(1,0)}M$ is a
holomorphic vector bundle with Hermitian metric induced by $g^{TM}$.
Let $\{e_{2i-1}\}_{i=1}^n\bigcup\{e_{2i}:=Je_{2i-1}\}_{i=1}^{n}$ be an
orthonormal frame of $TM$ and $\{e^i\}_{i=1}^{2n}$ be its dual
frame. Then
\begin{equation}\label{yy1}
\omega_j=\frac{1}{\sqrt{2}}(e_{2j-1}-\sqrt{-1}\,e_{2j}),\quad
j=1,\cdots,n
\end{equation}
form a local orthonormal frame of $T^{(1,0)}M$ with dual
frame$\{\omega^j\}_{j=1}^n$ (cf. \cite[(1.2.34)]{ma2007hmi}). We fix
this notation throughout the remaining part and use it without
further notice.

Let $\Theta$ be the real (1,1)-form defined by
\begin{equation}\label{4}
\Theta(u,v)=g^{TM}(Ju,v),\quad\text{for any\ } u,v\in \Gamma(TM).
\end{equation}
We call $\Theta$ as in {\rm (\ref{4})} a K\"{a}hler form on $M$. The
metric $g^{TM}$ on $TM$ is called a  K\"{a}hler metric and the
complex manifold $(M,J)$ is called a K\"{a}hler manifold if $\Theta$
is a closed form (cf. \cite[Definition 1.2.7]{ma2007hmi}).

There is a natural Hermitian metric on
$\Lambda^p\bigl(T^{*(1,0)}M\bigr)\otimes\Lambda^q\bigl(T^{*(0,1)}M\bigr)\otimes
E$ induced by $g^{TM}$ and $g^E$, which we denote by
$\langle\cdot,\cdot\rangle_{\Lambda^{*,*}\otimes E}$. By Wirtinger
Theorem (cf. \cite[pp. 31]{griffiths1994principles}), we know that
the volume form of $M$ determined by $g^{TM}$ is given by
$\tfrac{\Theta^n}{n!}$. Therefore, the $L^2$-scalar product on
$\Omega^{\ast,\ast}(M,E)$ is given by
\begin{equation}\label{6}
\langle\!\langle\alpha,\beta\rangle\!\rangle:=\int_{M}\langle\alpha,\beta\rangle_{
\Lambda^{*,*}\otimes E} \frac{\Theta^n}{n!}\ ,\qquad\text{for any
$\alpha,\beta\in\Omega^{\ast,\ast}(M,E)$.}
\end{equation}

The complex Hodge star operator is a complex conjugate linear
mapping $$\ast: \Omega^{p,q}(M,\C)\rightarrow\Omega^{n-p,
n-q}(M,\C)$$ such that if
$\alpha,\beta\in\Omega^{\ast,\ast}(M,\C)$, then (cf. \cite[pp.
80--82]{griffiths1994principles})
\begin{equation}
\alpha\wedge\ast\beta=\langle\alpha,\beta\rangle_{\Lambda^{*,*}}\frac{\Theta^n}{n!}.
\end{equation}
The formal adjoint $\overline{\partial}^\ast$ of
$\overline{\partial}$ with respect to the $L^2$-scalar product is
given by
\begin{equation}\label{8}
\overline{\partial}^\ast=-\ast\overline{\partial}\ast.
\end{equation}

Let $conj$ be the natural conjugate mapping induced by the bundle
automorphism (cf. \cite[pp. 141]{cappell2010cvant})
\begin{equation}
T^*M\otimes_\R\C\rightarrow T^*M\otimes_\R\C,\
v\otimes\lambda\mapsto v\otimes\overline{\lambda}\text{\ \ for any\
}v\in T^*M,\,\lambda\in\C.
\end{equation}
Then $\widehat{\ast}:=conj\,\ast$ is a complex linear mapping.
Clearly, $\widehat{\ast}=conj\,\ast=\ast\,conj$.

We now introduce the Clifford algebra (cf. \cite[Section
1.3.1]{ma2007hmi}).

For any $v\in TM$ with decomposition $v=v^{(1,0)}+v^{(0,1)}\in
T^{(1,0)}M\oplus T^{(0,1)}M$, let $\overline{v}^{(1,0),\ast}\in
T^{\ast(0,1)}M$ be the metric dual of $v^{(1,0)}$. The Clifford
action of $v$ on the bundle $\Lambda(T^{*(0,1)}M)=\Lambda^{\rm
even}(T^{*(0,1)}M)\oplus\Lambda^{\rm odd}(T^{*(0,1)}M)$ is defined
by
\begin{equation}
c(v)=\sqrt{2}\,\bigl(\overline{v}^{(1,0),\ast}\wedge-i_{v^{(0,1)}}\bigr),
\end{equation}
where $\wedge$ and $i$ denote the exterior and interior
multiplications, respectively. We verify easily that for $U,V\in
TM$,
\begin{equation}
c(U)c(V)+c(V)c(U)=-2\langle U,V\rangle.
\end{equation}
We identify $TM$ with $T^*M$ by the given metric $g^{TM}$, and
sometimes write $c(e^i)$ as $c(e_i)$. In the sequel, we don't
distinguish between $c(e^i)$ and $c(e_i)$.

Since $E$ is holomorphic, the usual $\overline{\partial}$ operator
on $\Omega^{\ast,\ast}(M,\C)$ has a unique natural extension to
$\Omega^{\ast,\ast}(M,E)$,
$\overline{\partial}_{E}:\Omega^{p,q}(M,E)\rightarrow\Omega^{p,q+1}(M,E)$
(cf. \cite[Section 3]{cappell2010cvant}).

Under the splitting
$\Omega^{1}(M,E)=\Omega^{1,0}(M,E)\oplus\Omega^{0,1}(M,E)$, the
connection $\nabla^E$ decomposes as sum:
$\nabla^E=(\nabla^E)^{1,0}\oplus(\nabla^E)^{0,1}$ with
\begin{equation}
(\nabla^E)^{1,0}: \Gamma(M,E)\rightarrow\Omega^{1,0}(M,E)\ ,\
(\nabla^E)^{0,1}: \Gamma(M,E)\rightarrow\Omega^{0,1}(M,E).
\end{equation}
Moreover, if we extend $\nabla^E$ on $\Gamma (M,E)$ in a unique way
to $\Omega^{\ast,\ast}(M,E)$ by Leibniz formula (cf. \cite[pp.
21]{berline2003hka}), then the extended $\nabla^E$ splits into two
pieces $\nabla^E=(\nabla^E)^{1,0}+(\nabla^E)^{0,1}$, which also
satisfy the Leibniz formula (cf. \cite[pp. 131]{berline2003hka}).
\begin{defn}{\rm (cf. \cite{cappell2010cvant})} The connection $\nabla^E$ is said to be compatible with the
holomorphic structure on $E$ if
$(\nabla^E)^{0,1}=\overline{\partial}_E$. The connection $\nabla^E$
is said to be of type $(1,1)$ if the curvature $(\nabla^E)^2$ is of
type $(1,1)$.
\end{defn}

\

In the sequel, we fix any $p$, $0\leq p\leq n$, and set
$\Omega^{p,\ast}(M,E)=\bigoplus_{q=0}^n\Omega^{p,q}(M,E)$.

Let $\nabla^E$ be a compatible (1,1) connection. Following
\cite[Section 3]{cappell2010cvant}, we define
\begin{equation}\label{24}
\overline{\partial}^{\ast}_{E,(\nabla^{E})^{1,0}}=-(\widehat{\ast}\otimes
1)(\nabla^{E})^{1,0}(\widehat{\ast}\otimes 1),
\end{equation}
and
\begin{equation}\label{25}
\Box_{E,\overline{\partial}}=\overline{\partial}_E\overline{\partial}^{\ast}_{E,(\nabla^{E})^{1,0}}
+\overline{\partial}^{\ast}_{E,(\nabla^{E})^{1,0}}\overline{\partial}_E.
\end{equation}
Since $\widehat{\ast}^2=\ast^2=(-1)^{p+q}$ on $\Omega^{p,q}(M,\C)$,
we get
\begin{equation}\label{28}
\left(\overline{\partial}^{\ast}_{E,(\nabla^{E})^{1,0}}\right)^2=(-1)^{p+q+1}(\widehat{\ast}\otimes
1)\left((\nabla^{E})^{1,0}\right)^2(\widehat{\ast}\otimes 1)=0.
\end{equation}

Using the following identities
\begin{equation}
\begin{split}
(\nabla^{E})^{1,0}&=\partial\otimes1+\omega^j\wedge\nabla^{E}_{\omega_j},\\
(\nabla^{E})^{0,1}&=\overline{\partial}\otimes1+\overline{\omega}^j\wedge\nabla^{E}_{\overline{\omega}_j},
\end{split}
\end{equation}
we get
\begin{equation}
\begin{split}\label{91}
\overline{\partial}^{\ast}_{E,(\nabla^{E})^{1,0}}&=\overline{\partial}^\ast\otimes1-i_{\overline{\omega}_j}\nabla^{E}_{\omega_j},\\
\overline{\partial}_E&=\overline{\partial}\otimes1+\overline{\omega}^j\wedge\nabla^{E}_{\overline{\omega}_j}.
\end{split}
\end{equation}

Set
\begin{equation}\label{45}
D=\sqrt{2}\left(\overline{\partial}_E+\overline{\partial}^{\ast}_{E,(\nabla^{E})^{1,0}}\right).
\end{equation}
Then from (\ref{25}), (\ref{91}), (\ref{45}), we deduce that
\begin{align}
\begin{split}\label{78}
D&=\sqrt{2}(\overline{\partial}+\overline{\partial}^\ast)\otimes1
+c(e_i)\otimes\nabla^{E}_{e_i},\\
D^2&=2\Box_{E,\overline{\partial}}\,.
\end{split}
\end{align}
\indent Denote ${\cal E}=\Lambda^p(T^{*(1,0)}M)\otimes E$. If we
assume in addition that $(M,J,g^{TM})$ is a K\"{a}hler manifold,
then using \cite[Proposition 3.67]{berline2003hka}, \cite[Lemma
1.4.4]{ma2007hmi}, the operator $D$ acting on
$\Omega^{p,\ast}(M,E)=\Omega^{0,\ast}(M,{\cal E})$ can be specified
as follows.
\begin{prop}
$D$ can be regarded as a Dirac operator on the Clifford module
$\Lambda(T^{*(0,1)}M)\otimes {\cal E}$ associated to the Clifford
connection
\begin{equation}\label{80}\nabla^{{\rm Cl}\otimes {\cal
E}}=\nabla^{{\rm Cl}}\otimes 1+1\otimes \nabla^{{\cal E}},\quad
\text{i.e.}\quad D=c(e_i)\nabla^{{\rm Cl}\otimes{\cal E}}_{e_i},
\end{equation}
where $\nabla^{{\rm Cl}}$ is the natural connection on
$\Lambda(T^{*(0,1)}M)$ induced by the Levi-Civita connection
$\nabla^{TM}$ {\rm (cf. \cite[pp. 28]{berline2003hka})} and
$\nabla^{\cal E}$ is the natural tensor connection induced by
$\nabla^{TM}$ and $\nabla^E$.
\end{prop}

\

Let ${\rm Spec}(\Box_{E,\overline{\partial}})$ be the spectrum of
$\Box_{E,\overline{\partial}}$ acting on $\Omega^{0,*}(M,{\cal E})$.
Denote
\begin{equation*}
{\rm Re}\bigl({\rm Spec}(\Box_{E,\overline{\partial}})\bigr)=\{{\rm
Re}(\lambda)\,|\,\lambda\in{\rm
Spec}(\Box_{E,\overline{\partial}})\}.
\end{equation*}
Pick any $a>0$ such that $a\notin{\rm Re}\left({\rm
Spec}(\Box_{E,\overline{\partial}})\right)$. Set
\begin{equation*}
\Omega_{<a}^{p,*}(M,E)=\bigoplus_{{\rm
Re}\,\lambda<a}\Omega_{\{\lambda\}}^{p,*}(M,E).
\end{equation*}
Let
$\Pi_{a}:\Omega^{p,*}(M,E)\longrightarrow\Omega_{<a}^{p,\ast}(M,E)$
be the orthogonal projection with respect to the natural
$L^2$-scalar product (cf. (\ref{6})) on $\Omega^{0,*}(M,{\cal E})$.
Set $P_{a}=1-\Pi_{a}$.

We recall some facts obtained in \cite[Section 4]{cappell2010cvant}.

$(\Omega_{<a}^{p,\ast}(M,E),\overline{\partial}_E)$ forms a finite
dimensional complex. Moreover, the inclusion of the complexes
\begin{equation*}
(\Omega_{<a}^{p,\ast}(M,E),\overline{\partial}_E)
\subset(\Omega^{p,\ast}(M,E),\overline{\partial}_E)
\end{equation*}
induces an isomorphism on cohomology. That is,
\begin{equation}\label{46}
H^q(\Omega_{<a}^{p,\ast}(M,E),\overline{\partial}_E)\stackrel{\cong}\longrightarrow
H^{p,q}_{\overline{\partial}}(M,E).
\end{equation}

$\widehat{\ast}\otimes1$ induces a $\C$-linear isomorphism of the
complex
$(\Omega_{<a}^{p,\ast}(M,E),\overline{\partial}^*_{E,(\nabla^E)^{1,0}})$
to the complex
$(\Omega_{<a}^{n-\ast,n-p}(M,E),(-1)^{p+1+*}(\nabla^E)^{1,0})$. We
have isomorphisms,
\begin{align}
\begin{split}\label{47}
H^{n-q,n-p}_{(\nabla^E)^{1,0}}(M,E)&\stackrel{\cong}\longrightarrow
H^{n-q}(\Omega_{<a}^{\ast,n-p}(M,E),(-1)^{p+1+*}(\nabla^E)^{1,0})\\&
\stackrel{\widehat{\ast}\otimes1}\longrightarrow
H_q\left(\Omega_{<a}^{p,\ast}(M,E),\overline{\partial}^*_{E,(\nabla^E)^{1,0}}\right).
\end{split}
\end{align}
\indent From \cite[Section 6]{cappell2010cvant}, we know that there
is a natural non-vanishing algebraic torsion invariant associated to
the complex
$\left(\Omega_{<a}^{p,\ast}(M,E),\overline{\partial}_E,\overline{\partial}^*_{E,(\nabla^E)^{1,0}}\right)$,
\begin{align*}
&{\rm
torsion}\left(\Omega_{<a}^{p,\ast}(M,E),\overline{\partial}_E,\overline{\partial}^*_{E,(\nabla^E)^{1,0}}\right)\\&\quad
\in\det\left(H^\ast\left(\Omega_{<a}^{p,\ast}(M,E),\overline{\partial}_E\right)\right)\otimes
\det\left(H_\ast\left(\Omega_{<a}^{p,\ast}(M,E),\overline{\partial}^*_{E,(\nabla^E)^{1,0}}\right)\right)^{-1}.
\end{align*}
Using the isomorphisms (\ref{46}) and (\ref{47}), we regard it as an
element of the complex line
\begin{equation*}
\det\left(H^{p,*}_{\overline{\partial}}(M,E)\right)\otimes
\det\left(H^{n-\ast,n-p}_{(\nabla^E)^{1,0}}(M,E)\right)^{-1}.
\end{equation*}
Note here the complex line is independent of the metric $g^{TM}$.\\
\indent Let $N$ be the number operator, i.e. $N$ acts on
$\Omega^{0,q}(M,\mathcal{E})$ by multiplication by $q$. By
\cite[Section 11]{cappell2010cvant}, for ${\rm Re}(s)>\tfrac{n}{2}$,
the following zeta-function is well-defined,
\begin{equation}
\zeta_{a}(s)=\tr_s\left[N(\Box_{E,\overline{\partial}})^{-s}P_{a}\right].
\end{equation}
Moreover, $\zeta_{a}(s)$ has a meromorphic extension to the whole
complex plane and is analytic at $s=0$. Consequently, the derivative
at $s=0$, $\zeta_{a}^\prime(0)$ is meaningful.

By \cite[Lemma 4.3]{cappell2010cvant} and
\cite[(4.8)]{cappell2010cvant}, we know that the combination
\begin{equation}\label{49}
{\rm
torsion}\left(\Omega_{<a}^{p,\ast}(M,E),\overline{\partial}_E,\overline{\partial}^*_{E,(\nabla^E)^{1,0}}\right)
\cdot\exp(\zeta_a^\prime(0))
\end{equation}
is independent of the choice of $a>0$ with $a\notin{\rm
Re}\left({\rm Spec}(\Box_{E,\overline{\partial}})\right)$ and the
Hermitian metric $g^E$.
\begin{defn}{\rm (cf. \cite[pp. 152]{cappell2010cvant})}
The non-vanishing element of the complex line
\begin{equation*}
\det\left(H^{p,*}_{\overline{\partial}}(M,E)\right)\otimes
\det\left(H^{n-\ast,n-p}_{(\nabla^E)^{1,0}}(M,E)\right)^{-1}
\end{equation*}
defined in {\rm (\ref{49})} is called the Cappell-Miller holomorphic
torsion and is denoted by $\tau_{{\rm holo},p}(M,E)$.
\end{defn}

\

%----------------------------------2.2-----------

\makeatletter
\@startsection{subsection}{0}{0pt}{2mm}{0.2mm}{\bfseries}{An anomaly
formula for the Cappell-Miller holomorphic torsion for
K\"{a}hler manifolds}\makeatother

We indicate the characteristic classes which we will use.

Let $g^{TM}$ be a K\"{a}hler metric on $TM$.  Let $R^+$ be the curvature of
$T^{(1,0)}M$ with the natural connection induced by the Levi-Civita
connection $\nabla^{TM}$ associated to $g^{TM}$. Let $R^E=(\nabla^E)^2$ be the curvature of
$\nabla^E$. If ${\cal A}\in \End(T^{(1,0)}M)$, set
\begin{equation}
\begin{split}
\td({\cal A})&={\det}_{T^{(1,0)}M}\Bigl(\frac{{\cal A}}{e^{{\cal
A}}-1}\Bigr),
\\
\det(I+t{\cal A})&=1+t\sigma_1({\cal A})+\dots +t^n{\sigma}_n({\cal
A}),
\\
\td_p({\cal A})&=\td({\cal A})\,\sigma_p(\exp{\cal A})\,.
\end{split}
\end{equation}

\begin{defn}Set
\begin{equation}
\td_p\bigl(T^{(1,0)}M,g^{TM}\bigr)=\td_p\left(\frac{R^+}{2\pi\sqrt{-1}}\right),\quad
\ch\bigl(E,\nabla^E\bigr)=\tr\left[\exp\left(\frac{-R^E}{2\pi\sqrt{-1}}\right)\right].
\end{equation}
\end{defn}

As in \cite[pp. 58]{bismut1988ata1}, let $P=\bigoplus_{j=0}^n\Omega^{j,j}(M,\C)$. Let $P'\subset P$ be the set
of smooth forms $\alpha\in P$ such that there exist smooth forms $\beta$, $\gamma$
for which $\alpha=\partial\beta+\overline{\partial}\gamma$. When $\alpha, \alpha'\in P$, we write $\alpha\equiv\alpha'$
if $\alpha-\alpha'\in P'$. Then the paring of the elements of $P/P'$ with the element of $P$ which is closed and has compact support
is well defined.

Let $g'^{TM}$ be another K\"{a}hler metric on $TM$.

By the results of \cite[Section (f)]{bismut1988ata1}, there is uniquely defined Bott-Chern class
\begin{equation*}
\widetilde{\td}_p\bigl(T^{(1,0)}M, g^{TM}, g'^{TM}\bigr)\in P/P'
\end{equation*}
such that
\begin{equation}
\frac{\overline{\partial}\partial}{2\pi\sqrt{-1}}\widetilde{\td}_p\bigl(T^{(1,0)}M, g^{TM}, g'^{TM}\bigr)
=\td_p\bigl(T^{(1,0)}M,g'^{TM}\bigr)-\td_p\bigl(T^{(1,0)}M,g^{TM}\bigr).
\end{equation}

Let $\tau_{{\rm holo},p}$, $\tau'_{{\rm holo},p}$ be the Cappell-Miller
holomorphic torsions associated to the K\"{a}hler metrics
$g^{TM}$, $g'^{TM}$, respectively.

Using the above notations, we state our main theorem as follows.

\begin{thm}\label{j2}{\rm (Compare with \cite[Theorem 1.23]{bismut1988ata3})}
The following identity holds,
\begin{equation}\label{j1}
\frac{\tau'_{{\rm holo},p}}{\tau_{{\rm holo},p}}=\exp\left(\,
\,\int_{M}\widetilde{\td}_p\bigl(T^{(1,0)}M, g^{TM}, g'^{TM}\bigr)
\cdot\ch\bigl(E,\nabla^E\bigr)\,\right)\,.
\end{equation}
\end{thm}

\

%--------------------------------3--------------------------------------------------------------------------------

\makeatletter
\@startsection{section}{0}{0pt}{2mm}{0.5mm}{\bfseries}{Some
estimates on heat kernels}
\makeatother

In this section, we deduce some estimates on heat kernels, which
will be needed in the next section.

Throughout this section, we assume that $M$ is a compact oriented smooth manifold of dimension $n$ with
a Riemannian metric $g^{TM}$ . Let $e_1,\cdots,e_n$ be a local
orthonormal frame of $TM$ with respect to $g^{TM}$. Let $C(M)$ be
the Clifford bundle associated to $g^{TM}$ (cf. \cite[Definition
3.30]{berline2003hka}). For any $0\leq i\leq n$, denote
\begin{equation}
\Omega^i(M)=\Gamma(\Lambda^i(T^*M)),\quad\text{and}\quad\Omega^*(M)=\bigoplus_{i=0}^n\Omega^i(M).
\end{equation}
Let ${\cal F}$ be a real vector bundle of dimension $m$ over $M$
with a connection $\nabla^{\cal F}$.

This section is self-contained.

%----------------------------3.1----------------

\makeatletter
\@startsection{subsection}{0}{0pt}{2mm}{0.2mm}{\bfseries}{The
special heat kernel depending on parameters} \makeatother

Let $\vartheta_1, \cdots, \vartheta_\imath$ be auxiliary Grassmann
variables, and any of which anticommutes with $C(X)$, where $C(X)$
is the element $X\in T_{x}M$ considered as an element of $C(M)$.
Assume also that the multiplication of any $q+1$ variables of the
above given Grassmann variables vanishes, where $q$ is some fixed
integer.

Let $R(\vartheta_1, \cdots, \vartheta_\imath)$ be the Grassmann
algebra generated by $1, \vartheta_1, \cdots, \vartheta_\imath$ (cf.
\cite{bismut1986aef2}). If $\omega\in R(\vartheta_1, \cdots,
\vartheta_\imath)$, then $\omega$ is a linear combinations of
$\vartheta_{i_1}\cdots\vartheta_{i_k}$, where $1\leq
i_1<\cdots<i_k\leq \imath$. We say the monomial
$\vartheta_{i_1}\cdots\vartheta_{i_k}$ is of degree $k$. Clearly,
$k\leq q$.

For any $t>0$, there is a homomorphism of algebras
\begin{equation}\label{23}
\psi_t: R(\vartheta_1, \cdots, \vartheta_\imath)\longrightarrow
R(\vartheta_1, \cdots, \vartheta_\imath),
\end{equation}
which for $1\leq j\leq \imath$, maps $\vartheta_j$ in
$\frac{\vartheta_j}{\sqrt t}$.

Defining the elements of $\Omega^*(M)\otimes\End({\cal F})$ to be of
degree zero, we give every monomial of $\Omega^*(M)\otimes\End({\cal
F})\widehat{\otimes}R(\vartheta_1, \cdots, \vartheta_\imath)$, say
$\varphi_{i_1\cdots i_k}\vartheta_{i_1}\cdots\vartheta_{i_k}$, where
$\varphi_{i_1\cdots i_k}\in\Omega^*(M)\otimes\End({\cal F})$, a
natural degree. The homomorphism in (\ref{23}) can be
extended to $\Omega^*(M)\otimes\End({\cal
F})\widehat{\otimes}R(\vartheta_1, \cdots, \vartheta_\imath)$ in an
obvious way.

For any $t>0$, $\omega\in\Omega^1(M)\otimes\End({\cal
F})\widehat{\otimes}R(\vartheta_1, \cdots, \vartheta_\imath)$, set
\begin{equation}
{}^t\omega=\psi_t(\omega),\qquad \text{and}\quad {}^t\nabla^{\cal
F}=\nabla^{\cal F}+{}^t\omega.
\end{equation}
We may assume that ${}^t\omega$ has no degree $0$ term \footnote{In
fact, let ${}^t\omega^{\{0\}}$ be the degree $0$ term of
${}^t\omega$, then $\widetilde{\nabla}^{\cal
F}:={}^t\omega^{\{0\}}+\nabla^{\cal F}$ is a connection on ${\cal
F}$. Replacing $\nabla^{\cal F}$ by $\widetilde{\nabla}^{\cal F}$,
then ${}^t\omega-{}^t\omega^{\{0\}}$ has no degree $0$ term.}.

For any $t>0$, $\rho\in\End({\cal F})\widehat{\otimes}R(\vartheta_1,
\cdots, \vartheta_\imath)$, set ${}^t\!\rho=\psi_t(\rho)$. We now
consider the smooth family of differential operators
\begin{equation*}
I_t^\prime=-\bigl(\nabla^{\cal
F}_{e_i}+{}^t\omega(e_i)\bigr)^2+{}^t\!\rho,\quad \text{and}\quad
I_t=tI_t^\prime.
\end{equation*}
Here we use the same notation as in \cite[Section
3(b)]{bismut1986asi},
\begin{equation}\label{21}
\bigl(\nabla^{\cal
F}_{e_i}+{}^t\omega(e_i)\bigr)^2:=\sum_{i=1}^{n}\bigl(\nabla^{\cal
F}_{e_i}+{}^t\omega(e_i)\bigr)^2-\nabla^{\cal
F}_{\nabla^{TM}_{e_i}e_i}-{}^t\omega(\nabla^{TM}_{e_i}e_i).
\end{equation}

For any $t>0$, let $k(x,y,s,t)$ be the smooth section of the bundle
${\cal F}\boxtimes{\cal F}^*$ over $\R_+\times M\times M$,
satisfying the following equation with boundary condition at $s=0$:
\begin{equation}\label{n1}
\left\{\begin{aligned}
&(\partial_s+I_t)k(x,y,s,t)=0,\\
&\lim_{s\rightarrow0+}\int_{y\in M}k(x,y,s,t)l(y) {\rm d}y=l(x).
\end{aligned}\right.
\end{equation}
Here we use the same convention as in \cite[pp. 72]{berline2003hka}
that $\mathcal {F}\boxtimes\mathcal {F}^*=\pi_1^*\mathcal
{F}\otimes\pi_2^*\mathcal {F}^*$, where $\pi_1$, $\pi_2$ are the
projections from $M\times M$ onto the first and second factor $M$
respectively.

Set $u=st$, then $\partial_s=t\partial_u$. The equation in
(\ref{n1}) is equivalent to
\begin{equation}\label{n2}
\left\{\begin{aligned}
&(\partial_u+I_t^\prime)\widetilde{k}(x,y,u,t)=0,\\
&\lim_{u\rightarrow0+}\int_{y\in M}\widetilde{k}(x,y,u,t)l(y){\rm
d}y=l(x),
\end{aligned}\right.
\end{equation}
where $\widetilde{k}(x,y,u,t)=k(x,y,s,t)$. In (\ref{n1}) and
(\ref{n2}), $l(x)$ is any smooth section of $\mathcal{F}$ and the
limit is meant in the uniform norm $\|l\|_0=\sup_{x\in M}\|l(x)\|$
for any metric on $\mathcal{F}$. $I_t^\prime$ is a generalized
Laplacian with parameters in $R(\vartheta_1,\cdots,
\vartheta_\imath)$ for any $t>0$.

Let $x$ and $y$ be sufficiently close points in $M$.  Let
$l_1,\cdots,l_{m}$ be a local frame of $\mathcal {F}$ on a
neighborhood of $y$, which are parallel along the radical geodesic
curve $x_s=\exp_{y}s{\bf x}:[0,1]\mapsto M$, with respect to the
connection $\nabla^{\cal F}$.

For any $t>0$, we define $\tau^t(x_s,y)\in {\rm Hom}({\cal
F}_y,{\cal F}_{x_s})$ as follows,
\begin{equation}
\tau^{t}(x_s,y)\bigl(l_1(y),\cdots,l_{m}(y)\bigr)=\bigl(l_1(x_s),\cdots,l_{m}(x_s)\bigr)A(s),
\end{equation}
where
\begin{equation}\label{1}
A(s)=I+\sum_{k=1}^q(-1)^k\int_{s\triangle_k}{}^t\omega(t_k)\cdots{}^t\omega(t_1){\rm
d}t_1\cdots{\rm d}t_k,
\end{equation}
and $s\triangle_k$, $k\geq1$, is a rescaled simplex given by
\begin{equation}
\{(t_1,\cdots t_k)|0\leq t_1\leq t_2 \cdots \leq t_k\leq s\}.
\end{equation}

Checking directly, we know that $A(s)$ is non-degenerated and that
\begin{equation}
{}^t\nabla^{\cal
F}_{\partial_s}\bigl(\tau^{t}(x_s,y)l_k(y)\bigr)=0,\quad \text{for
$k=1,\cdots,m$.}
\end{equation}
Therefore, it gives another smooth trivialization of $\cal F$ over
the neighborhood of $y$ corresponding to ${}^t\nabla^{\cal F}$.

 \indent Let $q_u$ be the smooth kernel modeled on
the Euclidean heat kernel: in a normal coordinates around $y$
$(x=\exp_y{\bf x})$,
\begin{equation*}
q_{u}(x,y)=(4\pi
u)^{-\frac{n}{2}}\exp\bigl(-\frac{\|\textbf{x}\|^2}{4u}\bigr).
\end{equation*}
We fix $y\in M$ and write $q_u$ for the section $x\mapsto q_u(x,y)$.

Set $j({\bf x})=\det^{\frac{1}{2}}(g_{ij}(\textbf{x}))$. Using
\cite[Theorem 2.26]{berline2003hka}, in the above new
trivialization, the formal solution\footnote{Because the manifold
considered here is oriented, there is a little but not essential
difference in the expression of the formal solution.} of the heat
equation (\ref{n2}) is given by
\begin{equation}
k_u(x,y;t)^{\rm
formal}=q_{u}(x,y)\sum_{i=0}^{\infty}u^i\Phi_i(x,y;t),
\end{equation}
where $\Phi_0(x,y;t)=j^{-\frac{1}{2}}(\textbf{x})\tau^t(x,y)$, and
for $i>0$,
\begin{equation*}
\Phi_i(x,y;t)=-j^{-\frac{1}{2}}(\textbf{x})\tau^t(x,y)\int_0^1s^{i-1}j^{\frac{1}{2}}(s
\textbf{x})\tau^t(x_s,y)^{-1}(I^\prime_{t,x}\Phi_{i-1})(x_s,y;t)\textrm{\rm
d}s.
\end{equation*}

Set
\begin{equation*}
\psi(s)=
\begin{cases}1,&\quad\text{if}\quad s<\frac{\varepsilon^2}{4},\\0,&\quad\text{if}\quad
s>\varepsilon^2,
\end{cases}
\end{equation*}
where $\varepsilon$ is chosen to be smaller than the injectivity
radius of the manifold $M$.

For $N$ large enough, we define the approximate solution by the
formula
\begin{equation}
k^N_u(x,y;t)=\psi\bigl(d(x,y)^2\bigr)q_u(x,y)\sum_{i=0}^Nu^i\Phi_i(x,y;t).
\end{equation}

In the sequel, we fix $N$ large enough. In a neighborhood of the
diagonal, we write $y=\exp_x \textbf{y}$, with $\textbf{y}\in T_xM$,
and identify $\mathcal {F}_x$ to $\mathcal {F}_y$ by parallel
transport along the geodesic joining $x$ and $y$ with respect to
$\nabla^{\cal F}$. That is, if $l$ is a smooth section of $\mathcal
{F}$, we denote by $l(x,\textbf{y})$ the function of $\textbf{y}$
such that $l(x,\textbf{y})=\tau(x,y)l(y)\in\mathcal {F}_x$; if
$\Phi$ is a section of $\mathcal {F}\boxtimes\mathcal {F}^*$, we
denote by $\Phi(x,\textbf{y})$ an endomorphism of $\mathcal {F}_x$
such that $\Phi(x,\textbf{y})=\Phi(x,y)\tau(x,y)^{-1}$. Let
$\Psi_i(x,\textbf{y};t)=\psi(\|\textbf{y}\|^2)\Phi_i(x,\textbf{y};t)$.
We have $\Psi_0(x,0)={\rm Id}_{\mathcal {F}_x}$.

The explicit formulas for $\Phi_i(x,y;t)$ and (\ref{1}) imply that
$\Psi_i(x,\textbf{y};t)$ can be expressed as
\begin{equation}\label{y13}
\Psi_i(x,{\bf y};t)=\sum_{j=0}^qt^{-\frac{j}{2}}\phi_{ij}(x,{\bf
y}),\quad\text{for all}\ i\geq 0,
\end{equation}
where $\phi_{ij}(x,{\bf y})\in\End({\cal F}_x)\widehat{\otimes}
R(\vartheta_1, \cdots, \vartheta_\imath)$.

Thus, the approximate solution can be expressed as
\begin{equation}
k_u^N(x,y;t)=q_{st}(x,y)\sum_{i=0}^N(st)^i\sum_{j=0}^qt^{-\frac{j}{2}}\phi_{ij}(x,{\bf
y})\tau(x,y).
\end{equation}

Let $(\mathscr{B}, \|\cdot\|)$ be a normed space. We introduced a norm
$\|\cdot\|_{\mathscr{B}\otimes R}$ on $\mathscr{B}\otimes R(\vartheta_1,
\cdots, \vartheta_\imath)$ as follows. For
\begin{equation*}
\varphi=\sum_{\tiny
\begin{subarray}{c} 1\leq k\leq q
\\
1\leq i_1<\cdots<i_k\leq \imath
\end{subarray}}
\varphi_{i_1\cdots i_k}\vartheta_{i_1}\cdots\vartheta_{i_k}\in \mathscr{B}
\otimes R(\vartheta_1, \cdots, \vartheta_\imath),\quad\text{we
define}
\end{equation*}
\begin{equation}\label{103}
 \|\varphi\|_{\mathscr{B}\otimes R}=\max_{\tiny
\begin{subarray}{c} 1\leq k\leq q
\\
1\leq i_1<\cdots<i_k\leq \imath
\end{subarray}
}\|\varphi_{i_1\cdots i_k}\|.
\end{equation}
Proceeding as \cite[Section 2.5]{berline2003hka}, we deduce that
there exists a unique smooth solution to (\ref{n1}), which is called
the heat kernel and is denoted by $p_u(x,y;t)$ (cf.
\cite[Proposition 2.17 and Theorem 2.30]{berline2003hka}).
Therefore, the unique smooth kernel $k(x,y,s,t)$ determined by
(\ref{n1}) equals to $p_{u}(x,y;t)|_{u=st}$. we can get the
following estimate.
\begin{thm}{\rm (compare with \cite[Theorem 2.23]{berline2003hka})}\label{y15}
Assume $0<s,t<T$. Then there exists a constant $C>0$, such that
\begin{equation}
\Bigl\|\partial_s^k\Bigl(p_u(x,y;t)-k_u^N(x,y;t)\Bigr)\Bigr\|_{\ell}\leq
C\cdot s^{N-\frac{n}{2}-\frac{\ell}{2}-k+1}\cdot
t^{N-{\frac{k+2}{2}}q-\frac{n}{2}-\frac{\ell}{2}+1}.
\end{equation}
where $\|\cdot\|_{\ell}$ is defined using the usual $\mathscr{C}^\ell$-norm {\rm (cf. \cite[pp. 71]{berline2003hka})} on $\mathscr{C}^\ell$-sections of ${\cal F}\boxtimes {\cal F}^*$ over
$\R_{+}\times M\times M$  as in {\rm (\ref{103})}.
\end{thm}

\

%---------------------------------3.2------------

\makeatletter
\@startsection{subsection}{0}{0pt}{2mm}{0.2mm}{\bfseries}{The
rescaled heat kernel and the Mehler formula} \makeatother

We continue the discussion of the last subsection.

Take any $x_0\in M$ and trivialize the vector bundle $\mathcal{F}$
in a neighborhood of $x_0$ by parallel transport along radical
geodesics with respect to $\nabla^{\cal F}$ (cf. \cite[pp.
153--154]{berline2003hka}). More precisely, let $V=T_{x_0}M$,
$F={\cal F}_{x_0}$ and $U=\{{\bf x}\in V|\|{\bf x}\|<\varepsilon\}$,
where $\varepsilon$ is the injectivity radius of the compact
manifold $M$ at $x_0$. We identify $U$ by means of exponential map
${\bf x}\mapsto\exp_{x_0}\!{\bf x}$ with a neighborhood of $x_0$ in
$M$. For $x=\exp_{x_0}\!{\bf x}$, the fibre ${\cal F}_{x}$ and $F$
are identified by the parallel transport map $\tau(x_0,x):{\cal
F}_x\rightarrow F$ along the geodesic $x_s=\exp_{x_0}\!s{\bf x}$.

Choose an orthonormal basis $\partial_i$ of $V=T_{x_0}M$, with dual
basis $d \textbf{x}^i$ of $T_{x_0}^*M$, and let $c^i=c(d
\textbf{x}^i)\in \End(E)$. Let $S$ be the spinor space of $V^*$ and
let $W=\Hom_{C(V^*)}(S,F)$ be the auxiliary vector space such that
$F=S\otimes W$, so that $\End(F)\cong\End(S)\otimes\End(W)\cong
C(V^*)\otimes \End(W)$. Let $e_i$ be the local orthonormal frame
obtained by parallel transports along the geodesics from the
orthonormal basis $\partial_i$ of $T_{x_0}M$, and let $e^i$ be the
dual frame of $T^*M$.

It follows that in this trivialization the bundle
$\End_{C(M)}\mathcal{F}$, restricted to $U$, is the trivial bundle
with fibre $\End_{C(V)}\mathcal{F}=\End(W)$ (cf. \cite[pp.
154]{berline2003hka}).

Let $p_u(x,x_0;t)$ be the smooth kernel of the operator $I_t$. We
transfer the kernel to the neighborhood $U$ of $0\in V$, thinking
of it as taking values in $\End(F)\widehat{\otimes} R(\vartheta_1,
\cdots, \vartheta_\imath)$, by writing
\begin{equation*}
p_u({\bf x};t)=\tau(x_0,x)p_u(x,x_0;t),\quad\text{where
$x=\exp_{x_0}\!{\bf x}$}.
\end{equation*}
Set
\begin{equation*}
\widehat{p}_u({\bf x};t)=\sigma(p_u({\bf x};t)),\quad\text{where
$\sigma$ is the full symbol map}.
\end{equation*}
Then $\widehat{p}_u({\bf x};t)$ is a $\Lambda(V^*)\otimes
\End(W)\widehat{\otimes} R(\vartheta_1, \cdots,
\vartheta_\imath)$-valued function on $U$. Consider the space
$\Lambda(V^*)\otimes \End(W)\widehat{\otimes} R(\vartheta_1, \cdots,
\vartheta_\imath)$ as a $C(V^*)\otimes \End(W)$ module, where the
action of $C(V^*)$ on $\Lambda(V^*)$ is the usual one
$c(e^j)=e^j\wedge-i_{e_j}$.

The $\Lambda V^*\otimes\End(W)\widehat{\otimes} R(\vartheta_1,
\cdots, \vartheta_\imath)$-valued function $\widehat{p}_u({\bf
x};t)$ satisfies the differential equation
\begin{equation}\label{y49}
(\partial_s+\widehat{I}_t({\bf x}))\widehat{p}_u({\bf x};t)=0,
\end{equation}
where $\widehat{I}_t$ is the local expression of $I_t$ in the above
trivialization.

Now we introduce Getzler rescaling (cf. \cite{getzler1986spl} and
\cite[Section 4.3]{berline2003hka}).

For any $\alpha\in\mathscr{C}^\infty\bigl(\mathbb{R}^+\times U,\Lambda
V^*\otimes\End_{C(V^*)}(F)\widehat{\otimes} R(\vartheta_1, \cdots,
\vartheta_\imath)\bigr)$, let
\begin{equation}\label{liu.rescaling}
(\delta_\varepsilon\alpha)(t,\textbf{x})=\sum_{i=0}^n\varepsilon^{-i}\alpha(\varepsilon^2t,\varepsilon\textbf{x})_{[i]},
\end{equation}
where $\alpha_{[i]}$ is the \emph{i}-form component of $\alpha$.

Set
$\widehat{I}_\varepsilon=\delta_\varepsilon\widehat{I}_t\delta_\varepsilon^{-1}$.
 Then we have
\begin{equation}\label{y36}
(\partial_s+\delta_\varepsilon\widehat{I}_t\delta_\varepsilon^{-1})(\delta_\varepsilon
\widehat{p}_u)=\delta_\varepsilon(\partial_s+\widehat{I}_t)\widehat{p}_u=0.
\end{equation}

\begin{defn}{\rm (compare with \cite[Definition 4.18]{berline2003hka})}
The rescaled heat kernel $r(\varepsilon,s,t,{\bf x})$ is defined by
\begin{equation}\label{y37}
r(\varepsilon,s,t,{\bf x})=\varepsilon^{n}(\delta_\varepsilon
\widehat{p}_u)(t,{\bf x}).
\end{equation}
\end{defn}

Fix any $T>1$. The following results are essentially similar to
those proved in \cite[Section 4.3]{berline2003hka}.
\begin{lem}{\rm (compare with \cite[Lemma 4.19]{berline2003hka})}\label{liu.rn}
There exist $\Lambda V^*\otimes\End(W)\widehat{\otimes} R(\vartheta_1, \cdots,
\vartheta_\imath)$\nobreakdash-valued polynomials $\gamma_i(s,t,{\bf
x})$ on $\mathbb{R}^+\times\mathbb{R^+}\times V$, such that for
every integer $N$, the function $r^N(\varepsilon,s,t,{\bf x})$
defined by
\begin{equation}
r^N(\varepsilon,s,t,{\bf x})=q_{st}({\bf
x})\sum_{i=-n-q}^{2N}\varepsilon^i\gamma_i(s,t,{\bf x})
\end{equation}
approximates $r(\varepsilon,s,t,{\bf x})$ in the following sense:

for $N$ large enough and $0<s,t<T$, there is a constant
$C(N,k,\alpha)>0$ such that
\begin{align*}
\left\|\partial_s^k\partial_{\bf
x}^\alpha\bigl(r(\varepsilon,s,t,{\bf x})-r^N(\varepsilon,s,t,{\bf
x})\bigr)\right\| <C(N,k,\alpha)\varepsilon^{2N-(k+2)q-n+2}\,,
\\
\text{for $(s,t,{\bf x})\in (0,T)\times(0,T)\times U$ and
$0<\varepsilon\leq 1$}\,,
\end{align*}
where $\|\cdot\|$ is defined by using the usual $\mathscr{C}^0$\nobreakdash-norm {\rm (cf. \cite[pp. 71]{berline2003hka})} on
$\mathscr{C}^0$\nobreakdash-sections of $\Lambda V^*\otimes\End(W)$
over $\R_{+}\times V$ as in {\rm (\ref{103})}.

Furthermore, for any $t\in (0,T)$, we have $\gamma_i(0,t,0)=0$ if
$i\neq0$, while $\gamma_0(0,t,0)={\rm Id}$.
\end{lem}

\

We state a general  Mehler formula for a generalized harmonic
oscillator.

Let $B$ be an $n\times n$ matrix, $L$ be an $m\times m$ matrix, both
with coefficients in the commutative algebra $\mathscr{A}$, where $n, m$
are any positive integers. The generalized harmonic oscillator is
the differential operator acting on $\mathscr{A}\otimes\End(\C^{m})$-valued functions defined by
\begin{equation*}
H=-\sum_i\bigl(\partial_i+\frac{1}{4}B_{ij}x_j\bigr)^2+L.
\end{equation*}

\begin{thm}\label{liu.mehler}{\rm (cf. \cite[Proposition 4.7]{donnelly1988lit})}
For any $a_0\in \End(\C^{m})$, there exists a unique formal solution
$p_u(x,B,L,a_0)$ of the heat equation
\begin{equation}\label{y35}
(\partial_u+H_x)p_u(x,B,L,a_0)=0
\end{equation}
of the form
\begin{equation}
p_u(x)=q_u(x)\sum_{k=0}^{\infty}u^k\Phi_k(x)
\end{equation}
and such that $\Phi_0(0)=a_0$. The function $p_u(x,B,L,a_0)$ is
given by the formula
\begin{equation*}
(4\pi u)^{-\frac{n}{2}}\widehat{\rm
A}(u{\cal D})\exp\Bigl(\frac{1}{8}{\cal C}_{ij}x_ix_j\Bigr)\exp\!\left(
-\frac{1}{4u}\Bigl(\frac{u{\cal D}}{2}\coth\frac{u{\cal D}}{2}
\Bigr)_{ij}x_ix_j\right)\exp(-uL)\,a_0,
\end{equation*}
where ${\cal C}_{ij}=\frac{1}{2}(B_{ij}+B_{ji})$,
${\cal D}_{ij}=\frac{1}{2}(B_{ij}-B_{ji})$ and ${\cal D}=({\cal D}_{ij})_{n\times n}$.
\end{thm}

\

%------------------------------------4----------------------------------------------------------------------------
\makeatletter \@startsection{section}{0}{0pt}{2mm}{0.5mm}{\bfseries}
{A proof of Theorem \ref{j2}}\makeatother

\makeatletter
\@startsection{subsection}{0}{0pt}{2mm}{0.2mm}{\bfseries}
{An infinitesimal variation formula for the Cappell-Miller
holomorphic torsion}\makeatother

Let $\ell\rightarrow g_\ell^{TM}$ be a smooth family of
K\"{a}hler metrics on $M$. Let $\ast_{\ell}$ be the complex Hodge star operators associated to the metrics
$g_{\ell}^{TM}$ acting on $\Omega^{p,\ast}(M)$. Let $D_{\ell}$ be
the operators acting on $\Omega^{p,\ast}(M,E)$ defined as in
(\ref{45}) corresponding to $ g_\ell^{TM}$ and $\nabla^E$.
Let $U_\ell=(g_{\ell}^{TM})^{-1}\frac{\partial}{\partial\ell}g_\ell^{TM}\in\End(TM\otimes_\mathbb{R}\mathbb{C})$.
Denote by $U_\ell^+\in {\rm End}(T^{(1,0)}M)$ the restriction of
$U_\ell$ to $T^{(1,0)}M$. Let $R^+_\ell$ be the curvature of
$T^{(1,0)}M$ with the natural connection induced by the Levi-Civita
connection $\nabla_\ell^{TM}$ associated to the K\"{a}hler metric
$g_{\ell}^{TM}$. Let $\tau_{{\rm holo},p,\ell}(M,E)$ be the
Cappell-Miller holomorphic torsion defined as in (\ref{49})
corresponding to $g_{\ell}^{TM}$.

The following infinitesimal variation formula for the Cappell-Miller
holomorphic torsion, which is closely related to \cite[Theorem
1.18]{bismut1988ata3}, is
obtained in \cite{cappell2010cvant}.
\begin{thm}\label{59}{\rm (cf. \cite[Lemma 4.2 and Lemma 7.1]{cappell2010cvant})}
As $t\rightarrow 0^+$, for every $k\in\mathbb{N}$, there is an
asymptotic expansion
\begin{equation}
\tr_s\Bigl[\ast_\ell^{-1}\frac{\partial{\ast_\ell}}{\partial
\ell}e^{-\frac{t}{2}{D^2_{\ell}}}\Bigr]=\sum_{j=-n}^{k}M_{j,{\ell}}\,t^j+o(t^k).
\end{equation}
Moreover,
\begin{equation}
\frac{\partial}{\partial\ell}{\rm log}\, \tau_{{\rm
holo},p,\ell}(M,E)=-M_{0,\ell}.
\end{equation}
\end{thm}

\

Let $\{e_j\}_{j=1}^{2n}$ be a local orthonormal frame of
$TM$ corresponding to $g^{TM}_\ell$, and
$\{\omega_{j}\}_{j=1}^{n}$ be the frame which is related to $\{e_j\}_{j=1}^{2n}$ as in
{\rm (\ref{yy1})}.

By direct calculation, we easily get the following lemma, of
which the $p=0$ case is proved in \cite[Proposition
1.19]{bismut1988ata3}.
\begin{lem}\label{79} The following identity holds,
\begin{equation}
\ast_{\ell}^{-1}\frac{\partial\ast_{\ell}}{\partial\ell}=\frac{\sqrt{-1}}{4}\dot{\Theta}(e_i,e_j)c(e_i)c(e_j)
+\frac{1}{4}\dot{\Theta}(e_i,Je_i)+\sqrt{-1}\
\dot{\Theta}(\omega_j,\overline{\omega}_i) \omega^j\wedge
i_{\omega_i}.
\end{equation}
\end{lem}

\

%-------------------------------------4.1--------------------------------
\makeatletter
\@startsection{subsection}{0}{0pt}{2mm}{0.2mm}{\bfseries} {The small
time asymptotics of the supertrace of certain heat
kernels}\makeatother

For simplicity, we set
\begin{equation}
Q_\ell=\ast_\ell^{-1}\frac{\partial{\ast_\ell}}{\partial \ell}.
\end{equation}

Take any $\ell$, say $\ell_0$. In the whole subsection, we always omit the subscript $\ell_0$ if there is no
confusion.

Then $Q$ splits into two parts,
\begin{align}
\begin{split}\label{63}
Q_{1}&=\frac{\sqrt{-1}}{4}\,\dot{\Theta}(e_i,e_j)c(e_i)c(e_j)
+\frac{1}{4}\,\dot{\Theta}(e_i,Je_i),\\
Q_{2}&=\sqrt{-1}\ \dot{\Theta}(\omega_j,\overline{\omega}_i)
\omega^j\wedge i_{\omega_i}.
\end{split}
\end{align}

\begin{prop}\label{81} As $t\rightarrow 0^+$, for every $k\in{\mathbb{N}}$, there are asymptotic
expansions
\begin{align}\label{101}
\begin{split}
\tr_s\Bigl[Q_{1}e^{-tD^2}\Bigr]=\sum_{j=-1}^{k}a_{j}\,t^j+o(t^k),\quad
\tr_s\Bigl[Q_{2}e^{-tD^2}\Bigr]=\sum_{j=0}^{k}b_{j}\,t^j+o(t^k).
\end{split}
\end{align}
Moreover, we have
\begin{align}
a_{-1}&=({2\pi
\sqrt{-1}})^{-n}\int_M\frac{\sqrt{-1}}{2}\,\dot{\Theta}\cdot\td(R^+)\cdot\tr\Bigl[\exp(-R^{\cal
E})\Bigr],
\\
b_{0}&=({2\pi \sqrt{-1}})^{-n}\int_M\sqrt{-1}\,\td(R^+) \cdot\tr\Bigl[\
\dot{\Theta}(\omega_i,\overline{\omega}_j) \omega^i\wedge
i_{\omega_j}\exp(-R^{\cal E})\Bigr],
\end{align}
where
\begin{equation}
R^{\cal E}=(\nabla^{\cal E})^2=-\langle
R^+\omega_i,\overline{\omega}_j\rangle_{g_{\ell_0}^{TM}}\omega^i\wedge
i_{w_j}+R^E.
\end{equation}
\end{prop}
\pf{Proof} The known result on the heat kernel asymptotic expansion
on the closed manifold $M$ (cf. \cite[Theorem 2.30]{berline2003hka})
implies the existence of the asymptotic expansion as
$t\rightarrow0$. On the other hand, using the standard local index
theory techniques of Getzler (cf. \cite{getzler1983pos},
\cite{getzler1986spl}, \cite[Section 4]{berline2003hka}), we deduce
from (\ref{63}) that
\begin{align*}
&\lim_{t\rightarrow 0+}\tr_s\Bigl[t\,Q_{1}e^{-tD^2}\Bigr]=({2\pi
i})^{-n}\int_M\frac{\sqrt{-1}}{2}\,\dot{\Theta}\cdot\td(R^+)\cdot\tr\Bigl[\exp(-R^{\cal
E})\Bigr],
\\
&\lim_{t\rightarrow 0+}\tr_s\Bigl[Q_{2}e^{-tD^2}\Bigr]=({2\pi
i})^{-n}\int_M\sqrt{-1}\,\td(R^+)
\tr\Bigl[\
\dot{\Theta}(\omega_i,\overline{\omega}_j) \omega^i\wedge
i_{\omega_j}\exp(-R^{\cal E})\Bigr].
\end{align*}
\indent We complete the proof of Proposition \ref{81}.

\

The remaining part of this subsection is devoted to  the calculation of
$a_0$. We follow \cite[Section 1(h)]{bismut1988ata3} in spirt,
although we will not use probability theory in our final
calculation. Instead, using the estimates of heat kernels with
parameters deduced in last section, we modify the proof of
the local index theorem presented in \cite[Chapter 4]{berline2003hka}.

We use the convention as in \cite[Section
2(f)]{bismut1988ata2}. Let ${\rm d}a$, ${\rm d}\overline{a}$ be two
odd Grassmann variables. If $\eta\in
\Lambda\bigl(T^*M\otimes_{\mathbb{R}}\mathbb{C}\bigr)\widehat{\otimes}\C({\rm d}a,{\rm d}\overline{a})$
can be written in the form
\begin{equation*}
\eta=\eta_0+{\rm d}a\,\eta_1+{\rm d}\overline{a}\,\eta_2+{\rm
d}a{\rm d}\overline{a}\,\eta_3,\quad\text{where}\ \eta_i\in
\Lambda\bigl(T^*M\otimes_{\mathbb{R}}\mathbb{C}\bigr),0\leq i\leq 3,
\end{equation*}
then we set
\begin{equation}(\eta)^{{\rm d}a{\rm
d}\overline{a}}=\eta_3.
\end{equation}

The following two propositions are proved in \cite{bismut1988ata3}.
\begin{prop}{\rm (cf. \cite[Theorem 1.20]{bismut1988ata3})}\label{75} The
following identity holds,
\begin{equation}\label{remark}
\begin{split}
&\ \frac{\partial}{\partial
t}\Bigl(t\tr_s\Bigl[Q_{1}\exp(-tD^2)\Bigr]\Bigr)\\
&=\Bigl(\tr_s\Bigl[\exp\Bigl(-tD^2-\sqrt{\frac{t}{2}}{\rm
d}aD-\sqrt{\frac{t}{2}}{\rm d}\overline{a}[D,Q_{1}]+{\rm d}a{\rm
d}\overline{a}\,Q_{1}\Bigr)\Bigr]\Bigr)^{{\rm d}a{\rm
d}\overline{a}}.
\end{split}
\end{equation}
\end{prop}

\begin{prop}{\rm(cf. \cite[Theorem 1.21]{bismut1988ata3})}\label{74} The following identity holds,
\begin{equation}
\begin{split}\label{yb1}
&-tD^2-\sqrt{\frac{t}{2}}\,{\rm d}aD-\sqrt{\frac{t}{2}}\,{\rm
d}\overline{a}[D,Q_{1}]+{\rm d}a{\rm d}\overline{a}\,Q_{1}
\\
&=t\Bigl(\nabla^{{\rm Cl}\otimes {\cal
E}}_{e_i}-\frac{1}{2\sqrt{2t}}{\rm d}ac(e_i)-\frac{{\rm
d}\overline{a}}{2\sqrt{2t}}\sqrt{-1}\,\dot\Theta(e_k,e_i)c(e_k)\Bigr)^2\\
&\qquad+\frac{{\rm d}a{\rm
d}\overline{a}}{4}\dot\Theta(e_j,Je_j)-\sqrt{\frac{t}{2}}{\rm
d}\overline{a}\frac{c(e_i)}{4}(\nabla^{TM}_{e_i}\dot\Theta)(e_j,Je_j)\\
&\ \qquad-\frac{tK}{4}-\frac{t}{2}c(e_i)c(e_j)F^{{\cal
E}\!/\!S}(e_i,e_j),
\end{split}
\end{equation}
where $K$ is the scalar curvature of $M$, and $F^{{\cal E}\!/\!S}$
denotes the twisting curvature of the Clifford module
$\Lambda\bigl({T^{*(0,1)}M}\bigr)\otimes{\cal E}$ given by {\rm (cf.
\cite[pp. 117, 148]{berline2003hka})}
\begin{equation}\label{twist}
F^{\mathcal{E}/\!S}=\frac{1}{2}\tr_{T^{(1,0)}M}[R^+]+R^{\cal E}.
\end{equation}
In {\rm (\ref{yb1})} we also use
the same notations as in {\rm \cite[Section 3(b)]{bismut1986asi}
(compare with (\ref{21}))}.
\end{prop}

\

Set
\begin{equation}\label{82}
I_t=tD^2+\sqrt{\frac{t}{2}}\,{\rm d}aD+\sqrt{\frac{t}{2}}\,{\rm
d}\overline{a}[D,Q_{1}]-{\rm d}a{\rm d}\overline{a}\,Q_{1}.
\end{equation} For the bundle
$\Lambda\bigl({T^{*(0,1)}M}\bigr)\otimes{\cal E}$, by Proposition
\ref{74}, we  see that $I_t$ is in the type of operators we study in
last section, and all the results we get there can be applied to
the current case.

\begin{lem}\label{n3} The following identity holds,
\begin{equation}
\tr_s[\exp(-I_t)]=\int_M\tr_s[k(x,x,1,t)]{\rm d}v_M(x),
\end{equation}
where $k(x,y,s,t)$ is the unique solution of {\rm (\ref{n1})}.
\end{lem}

\

We trivialize the bundle
$\Lambda\bigl({T^{*(0,1)}M}\bigr)\otimes{\cal E}$ by parallel
transports along the radical geodesics with respect to the Clifford
connection $\nabla^{{\rm Cl}\otimes{\cal E}}$ as in last section. Using
\cite[Lemma 4.14 and Lemma 4.15]{berline2003hka} and Proposition
\ref{74}, we can easily deduce the local expression of $I_t$.
\begin{lem}
In the chosen trivialization, the operator $I_t$, when restricted to
$U$, can be specified as the following operator $\widehat{I}_t$,
which is a differential operator on $U$ with coefficients in
$C(V^*)\otimes\End(W)\otimes R({\rm d}a, {\rm d}\overline{a})$,
\begin{align}
\widehat{I}_t&=-tg^{ij}({\bf x})\left(\nabla^{{\rm Cl}\otimes{\cal
E}}_{\partial_i}-\frac{{\rm
d}a}{2\sqrt{2t}}\langle\partial_i,e_k\rangle({\bf x})c^k
-\frac{\sqrt{-1}}{2\sqrt{2t}}{\rm
d}\overline{a}\cdot\dot{\Theta}(e_k,\partial_i)({\bf
x})c^k\right)\notag\\
&\qquad\qquad\cdot\left(\nabla^{{\rm Cl}\otimes{\cal
E}}_{\partial_j}-\frac{{\rm
d}a}{2\sqrt{2t}}\langle\partial_j,e_l\rangle({\bf x})c^l
-\frac{\sqrt{-1}}{2\sqrt{2t}}{\rm
d}\overline{a}\cdot\dot{\Theta}(e_l,\partial_j)({\bf
x})c^l\right)\notag\\
&\quad+tg^{ij}({\bf x})\Gamma_{ij}^k({\bf x})\left(\nabla^{{\rm
Cl}\otimes{\cal E}}_{\partial_k}-\frac{{\rm
d}a}{2\sqrt{2t}}\langle\partial_k,e_l\rangle({\bf x})c^l
-\frac{\sqrt{-1}}{2\sqrt{2t}}{\rm
d}\overline{a}\cdot\dot{\Theta}(e_l,\partial_k)({\bf
x})c^l\right)\notag\\&\quad-\frac{1}{4}{\rm d}a{\rm
d}\overline{a}\cdot\dot{\Theta}(e_i,Je_i)({\bf
x})+\frac{1}{4}\sqrt{\frac{t}{2}}{\rm
d}\overline{a}\cdot(\nabla^{TM}_{e_i}\dot{\Theta})(e_j,Je_j)({\bf
x})c^i \notag\\&\quad+\frac{tK}{4}({\bf
x})+\frac{t}{2}F^{\mathcal{E}\!/\!S}(e_i,e_j)({\bf x})c^ic^j.
\end{align}
\end{lem}

We compute the Getzler rescaling (cf. (\ref{liu.rescaling})) of the
operator $\widehat{I}_t$.
\begin{align}
\widehat{I}_\varepsilon=&-\varepsilon^2t
g^{ij}(\varepsilon\textbf{x})\left(\delta_\varepsilon\nabla^{{\rm
Cl}\otimes{\cal
E}}_{\partial_i}\delta_\varepsilon^{-1}-\frac{\varepsilon^{-1}}{2\sqrt{2t}}{\rm
d}a\cdot \langle\partial_i,e_k\rangle(\varepsilon\textbf{x})\cdot
(\varepsilon^{-1}e^k\wedge-\varepsilon i_{e_k})\right.\notag\\&
\hspace{90pt}\left.-\frac{\sqrt{-1}\varepsilon^{-1}}{2\sqrt{2t}}{\rm
d}\overline{a}\cdot
\dot{\Theta}(e_k,\partial_i)(\varepsilon\textbf{x})\cdot
(\varepsilon^{-1}e^k\wedge-\varepsilon
i_{e_k})\right)\notag\\
&\hspace{30pt}\cdot\left(\delta_\varepsilon\nabla^{{\rm
Cl}\otimes{\cal
E}}_{\partial_j}\delta_\varepsilon^{-1}-\frac{\varepsilon^{-1}}{2\sqrt{2t}}{\rm
d}a\cdot \langle\partial_j,e_l\rangle(\varepsilon\textbf{x})\cdot
(\varepsilon^{-1}e^l\wedge-\varepsilon i_{e_l})\right.\notag\\
&\hspace{50pt}\left.-\frac{\sqrt{-1}\varepsilon^{-1}}{2\sqrt{2t}}{\rm
d}\overline{a}\cdot
\dot{\Theta}(e_l,\partial_j)(\varepsilon\textbf{x})\cdot
(\varepsilon^{-1}e^l\wedge-\varepsilon i_{e_l})\right)\notag\\
&+\varepsilon^2tg^{ij}(\varepsilon\textbf{x})\Gamma_{ij}^k(\varepsilon\textbf{x})
\left(\delta_\varepsilon\nabla^{{\rm Cl}\otimes{\cal
E}}_{\partial_k}\delta_\varepsilon^{-1}-\frac{\varepsilon^{-1}}{2\sqrt{2t}}{\rm
d}a\cdot
\langle\partial_k,e_l\rangle(\varepsilon\textbf{x})(\varepsilon^{-1}e^l\wedge-\varepsilon i_{e_l})\right.\notag\\
&\left.\hspace{40pt}
-\frac{\sqrt{-1}\varepsilon^{-1}}{2\sqrt{2t}}{\rm
d}\overline{a}\cdot
\dot{\Theta}(e_l,\partial_k)(\varepsilon\textbf{x})\cdot
(\varepsilon^{-1}e^l\wedge-\varepsilon i_{e_l})\right)\notag\\
&-\frac{1}{4}{\rm d}a{\rm d}\overline{a}\!\cdot\!
\dot{\Theta}(e_i,Je_i)(\varepsilon\textbf{x})
+\frac{\varepsilon}{4}\sqrt{\frac{t}{2}}{\rm
d}\overline{a}\!\cdot\!(\nabla^{TM}_{e_i}\dot{\Theta})(e_j,Je_j)
(\varepsilon\textbf{x})\!\cdot\!(\varepsilon^{-1}e^i\wedge-\varepsilon
i_{e_i})\notag\\
&+\frac{\varepsilon^2tK}{4}(\varepsilon\textbf{x})
+\frac{\varepsilon^2t}{2}F^{\mathcal{E}\!/\!S}(e_i,e_j)(\varepsilon\textbf{x})
(\varepsilon^{-1}e^i\wedge-\varepsilon
i_{e_i})(\varepsilon^{-1}e^j\wedge-\varepsilon i_{e_j}),
\end{align}
where (cf. \cite[Lemma 4.15]{berline2003hka})
\begin{align}
\delta_\varepsilon\nabla^{{\rm Cl}\otimes{\cal
E}}_{\partial_i}\delta_\varepsilon^{-1}&=\varepsilon^{-1}\partial_i+
\frac{\varepsilon}{8}R_{klij}\textbf{x}^j(\varepsilon^{-1}e^k\wedge-\varepsilon
i_{e_k}) (\varepsilon^{-1}e^l\wedge-\varepsilon i_{e_l})\notag
\\
&+\frac{1}{2}f_{ikl}(\varepsilon\textbf{x})(\varepsilon^{-1}e^k\wedge-\varepsilon
i_{e_k}) (\varepsilon^{-1}e^l\wedge-\varepsilon
i_{e_l})+g_i(\varepsilon\textbf{x}),
\\
R_{klij}=\bigl(R&(\partial_i,\partial_j)\partial_l,\partial_k\bigr)_{x_0}\text{
is the Riemannian curvature at $x_0$, and}\notag
\\
f_{ikl}(\textbf{x})&=O(|\textbf{x}|^2)\in C^{\infty}(U),\ g_i({\bf
x})=O(|{\bf x}|)\in \Gamma\bigl(U,\End(W)\bigr). \notag
\end{align}

The singular term of $\widehat{I}_\varepsilon$, as
$\varepsilon\rightarrow0$, is that
\begin{equation}
\begin{split}
&t\varepsilon^{-1}\Bigl(\partial_i+R_{klij}\textbf{x}^je^k\wedge
e^l\wedge+\frac{1}{2}\varepsilon^{-2}f_{ikl}(\varepsilon\textbf{x})e^k\wedge
e^l\wedge\Bigr)
\\
&\qquad\cdot\Bigl(\frac{{\rm
d}a}{\sqrt{2\,t}}\cdot\langle\partial_i,e_k\rangle(0)
e^k\wedge+\frac{\sqrt{-1}}{\sqrt{2\,t}}{\rm
d}\overline{a}\cdot\dot{\Theta}(e_k,\partial_i)(0)e^k\wedge\Bigr)
\\
&+t\varepsilon^{-2}\Bigl(\frac{\sqrt{-1}}{4t}{\rm d}a{\rm
d}\overline{a}\cdot\langle\partial_i,e_k\rangle(0)\cdot\dot{\Theta}(e_l,\partial_i)(0)e^k\wedge
e^l\wedge\Bigr).
\end{split}
\end{equation}
To remove it, we proceed as in \cite{donnelly1988lit} by
conjugation. Take
\begin{equation*}
A(\textbf{x},\varepsilon)=\sum_{i,k}\left(\frac{\varepsilon^{-1}}{2\sqrt{2t}}{\rm
d}a
\langle\partial_i,e_k\rangle(0)\textbf{x}^ie^k\wedge+\frac{\sqrt{-1}}
{2\sqrt{2t}}\varepsilon^{-1}{\rm
d}\overline{a}\cdot\dot{\Theta}(e_k,\partial_i)(0)\textbf{x}^ie^k\wedge\right),
\end{equation*}
and set $h=\exp(-A)$. Clearly, $h$ is a polynomial of $\textbf{x}$ and $h(0)=1$.\\
\indent Set $$\widehat{J}_\varepsilon=h\widehat{I}_\varepsilon
h^{-1}.$$ From (\ref{y36}) and (\ref{y37}), we get
\begin{equation}\label{y38}
(\partial_s+\widehat{J}_\varepsilon)\Bigl(h({\bf x})\cdot
r(\varepsilon,s,t,{\bf x})\Bigr)=h({\bf
x})\cdot(\partial_s+\widehat{I}_\varepsilon)r(\varepsilon,s,t,{\bf
x})=0.
\end{equation}
Calculating directly, we deduce
\begin{align}
\widehat{J}_\varepsilon=&-t\left(\partial_i+\frac{1}{8}R_{klij}\textbf{x}^je^k\wedge
e^l-\frac{\varepsilon^{-1}}{2\sqrt{2t}}{\rm d}a
\Bigl(\langle\partial_i,e_k\rangle(\varepsilon\textbf{x})-\langle\partial_i,e_k\rangle(0)\Bigr)e^k\wedge\right.\notag\\
&\hspace{40pt}-\frac{\sqrt{-1}\varepsilon^{-1}}{2\sqrt{2t}}{\rm
d}\overline{a}\cdot
\left(\dot{\Theta}(e_k,\partial_i)(\varepsilon\textbf{x})-\dot{\Theta}(e_k,\partial_i)(0)\right)e^k\wedge\notag\\
&\hspace{60pt}\left.-\frac{\sqrt{-1}}{4t}{\rm d}a{\rm
d}\overline{a}\cdot\dot{\Theta}(\partial_i,\partial_j)(0)\textbf{x}^j\right)^2\notag\\
&-\frac{1}{4}{\rm d}a{\rm
d}\overline{a}\cdot\dot{\Theta}(e_i,Je_i)(0)
+\frac{1}{4}\sqrt{\frac{t}{2}}{\rm
d}\overline{a}\cdot(\nabla^{TM}_{e_i}\dot{\Theta})(e_j,Je_j)
(0)e^i\wedge\notag\\
&+\frac{t}{2}F^{\mathcal{E}\!/\!S}(e_i,e_j)(0)e^i\wedge e^j+
error(\textbf{x},\varepsilon).
\end{align}
The symbol $error(\textbf{x},\varepsilon)$ denotes the terms which
vanish when $\varepsilon\rightarrow0$ and will not contribute in the
final analysis.\\
\indent Clearly, for fixed $\textbf{x}\in U$,
$\lim_{\varepsilon\rightarrow0^+}\widehat{J}_\varepsilon$ exists and
is given by
\begin{equation}
\widehat{J}_0=-t\sum_i\Bigl(\partial_i+\frac{1}{4}B_{ij}\textbf{x}^j\Bigr)^2+tL,
\end{equation}
where
\begin{align}\label{y44}
B_{ij}=&\frac{1}{2}R_{klij}e^k\wedge e^l-
\sqrt{\frac{2}{t}}\sqrt{-1}{\rm d}\overline{a}\cdot
\bigl(\nabla^{TM}_{\partial_j}\dot{\Theta}\bigr)(e_k,\partial_i)(0)e^k\wedge\notag\\
&-\frac{\sqrt{-1}}{t}{\rm d}a{\rm
d}\overline{a}\cdot\dot{\Theta}(\partial_i,\partial_j)(0),\\
\label{y45} tL=&-\frac{1}{4}{\rm d}a{\rm
d}\overline{a}\cdot\dot{\Theta}(e_i,Je_i)(0)
+\frac{1}{4}\sqrt{\frac{t}{2}}{\rm
d}\overline{a}\cdot(\nabla^{TM}_{e_i}\dot{\Theta})(e_j,Je_j)
(0)e^i\wedge\notag\\
&+\frac{t}{2}F^{\mathcal{E}\!/\!S}(e_i,e_j)(0)e^i\wedge e^j.
\end{align}
In the above calculation of
$\lim_{\varepsilon\rightarrow0^+}\widehat{J}_\varepsilon$, we use
\cite[Proposition 1.28]{berline2003hka}, which claims that as
$\varepsilon\rightarrow 0^+$, $\langle\partial_i,e_k\rangle
(\varepsilon\textbf{x})=\delta_{ik}+O(\varepsilon^2)$. Therefore, the term
containing ${\rm d}a$ in $\widehat{J}_\varepsilon$ vanishes as
$\varepsilon\rightarrow 0^+$.

Using the fact that
$\widehat{J}_{\varepsilon}=\widehat{J}_0+O(\varepsilon)$, we can now
show that there are no poles in the Laurent series expansion in
$\varepsilon$ of  $h({\bf x})r(\varepsilon,s,t,{\bf x})$. By Lemma
\ref{liu.rn}, we have
\begin{equation}
h({\bf x})r(\varepsilon,s,t,{\bf x})\sim
q_{st}(\textbf{x})\sum_{i=-2n-2}^{\infty}\varepsilon^ih({\bf
x})\gamma_i(s,t,\textbf{x}).
\end{equation}
\indent We expand the equation (\ref{y38})
\begin{equation*}
(\partial_s+\widehat{J}_\varepsilon)\Bigl(h({\bf x})\cdot
r(\varepsilon,s,t,{\bf x})\Bigr)=0
\end{equation*}
in a Laurent series in $\varepsilon$. Lemma \ref{liu.rn} implies
that the leading term
$$q_{st}(\textbf{x})\varepsilon^{-l}h({\bf
x})\gamma_{-l}(s,t,\textbf{x})$$ of the asymptotic expansion of
$h({\bf x})r(\varepsilon,s,t,{\bf x})$ satisfies the heat equation
\begin{equation}
(\partial_s+\widehat{J}_0)\Bigl(q_{st}(\textbf{x})h({\bf
x})\gamma_{-l}(s,t,\textbf{x})\Bigr)=0,\ \text{for any fixed small}\
t>0,
\end{equation}
which is equivalent to
\begin{equation}\label{y39}
\Bigl(\partial_u+\frac{1}{t}\widehat{J}_0\Bigl)\Bigl(q_{st}(\textbf{x})h({\bf
x})\gamma_{-l}(s,t,\textbf{x})\Bigr)=0,\ \text{where $u=st$}.
\end{equation}
Since $\frac{1}{t}\widehat{J}_0$ is a harmonic oscillator, we can
apply the generalized Mehler formula Theorem \ref{liu.mehler} to
(\ref{y39}). The boundary condition $\gamma_{-l}(0,t,0)=0$ for $l>0$
implies $\gamma_{-l}(s,t,\textbf{x})\equiv0$ for $l>0$. In
particular, we see that there are no poles in the Laurent series
expansion of $h({\bf x})r(\varepsilon,s,t,{\bf x})$ in powers of
$\varepsilon$.\\
\indent The other thing that we learn from the above argument is
that the leading term of the expansion of $h({\bf
x})r(\varepsilon,s,t,{\bf x})$, i.e. $h({\bf x})r(0,s,t,{\bf
x})=q_{st}(\textbf{x})h({\bf x})\gamma_0(s,t,\textbf{x})$, satisfies
the equation
\begin{equation}
\Bigl(\partial_u+\frac{1}{t}\widehat{J}_0\Bigl)\Bigl(q_{st}(\textbf{x})h({\bf
x})\gamma_{0}(s,t,\textbf{x})\Bigr)=0,\ \text{with}\
\gamma_{0}(0,t,0)=1.
\end{equation}
Using Theorem \ref{liu.mehler},  $q_{st}(\textbf{x})h({\bf
x})\gamma_{0}(s,t,\textbf{x})$ is given by the explicit formula
\begin{equation}
(4\pi u)^{-n}\widehat{\rm A}(u{\cal D})\exp\Bigl(\frac{1}{8}\,{\cal C}_{ij}{\bf
x}^i{\bf x}^j\Bigr)\exp\left(-uL
-\frac{1}{4u}\Bigl(\frac{u{\cal D}}{2}\coth\frac{u{\cal D}}{2} \Bigr)_{ij}{\bf
x}^i{\bf x}^j\right)\notag.
\end{equation}
In particular, we get for any fixed small $t>0$ that
\begin{equation}\label{y40}
\lim_{\varepsilon\rightarrow0+}r(\varepsilon,1,t,0)=r(0,1,t,0)=(4\pi
t)^{-n}\widehat{\rm A}(t{\cal D})\exp(-tL).
\end{equation}
\indent By (\ref{liu.rescaling}) and (\ref{y37}), we have
\begin{equation*}
r(\varepsilon,s,t,{\bf
x})=\sum_{i=0}^{2n}\varepsilon^{2n-i}\widehat{p}_{\varepsilon^2ts}(\varepsilon{\bf
x};\varepsilon^2t)_{[i]},
\end{equation*}
from which, we obtain
\begin{equation}
\begin{split}\label{y41}
r(\varepsilon,s,t,{\bf x})|_{s=1,{\bf
x}=0}&=\sum_{i=0}^{2n}\varepsilon^{2n-i}\widehat{p}_{\varepsilon^2t}(0;\varepsilon^2t)_{[i]},\\
\Bigl(r(\varepsilon,1,t,0)\Bigr)_{[2n]}&=\widehat{p}_{\varepsilon^2t}(0;\varepsilon^2t)_{[2n]}.
\end{split}
\end{equation}
Moreover,
\begin{equation*}
\lim_{\varepsilon\rightarrow0^+}
\tr\bigl[\widehat{p}_{\varepsilon^2t}(0;\varepsilon^2t)_{[2n]}\bigr]\
\text{exists}.
\end{equation*}

\indent Using \cite[Proposition 3.21]{berline2003hka}, we get for
any fixed small $t>0$ that
\begin{equation}\label{y42}
\tr_s\bigl[k(x_0,x_0,1,t)\bigr]{\rm d}v_M(x_0)
=(-2\sqrt{-1})^{n}\tr\bigl[\widehat{p}_t(0;t)_{[2n]}\bigr].
\end{equation}
\indent Since
\begin{equation*}
\lim_{t\rightarrow0^+}\tr\,\Bigl[\widehat{p}_t(0;t)_{[2n]}\Bigr]=\lim_{\varepsilon\rightarrow0^+}
\tr\,\Bigl[\widehat{p}_{\varepsilon^2t}(0;\varepsilon^2t)_{[2n]}\Bigr],
\end{equation*} from (\ref{y40}), (\ref{y41}) and (\ref{y42}), we
get
\begin{equation}
\lim_{t\rightarrow0+}\tr_s\bigl[k(x_0,x_0,1,t)\bigr]{\rm d}v_M(x_0)=(2\pi
\sqrt{-1}\,t)^{-n}\tr\Bigl[\widehat{\rm A}(t{\cal D})\exp(-tL)\Bigr]_{[2n]}.
\end{equation}
Moreover, by Lemma \ref{n3}, we deduce
\begin{equation}\label{y43}
\lim_{t\rightarrow0}\tr_s\left[e^{-I_t}\right]=(2\pi
\sqrt{-1}\,t)^{-\frac{n}{2}}\int_M\tr\Bigl[\widehat{\rm A}(t{\cal D})\exp(-tL)\Bigr].
\end{equation}

From (\ref{twist}), (\ref{y44}), (\ref{y45}) and (\ref{y43}), we have
\begin{align}
\begin{split}
&\lim_{t\rightarrow0}\left(\tr_s[e^{-I_t}]\right)^{{\rm d}a{\rm
d}\overline{a}}=(2\pi \sqrt{-1}\,t)^{-n}\left(\int_M\tr\Bigl[\widehat{\rm
A}(t{\cal D})\exp(-tL)\Bigr]\right)^{{\rm d}a{\rm d}\overline{a}}
\\
&=(2\pi\sqrt{-1})^{-n}\left(\int_M\tr\left[\widehat{\rm
A}\Bigl(\frac{1}{2}R_{klij}e^k\wedge e^l -\sqrt{-1}{\rm d}a{\rm
d}\overline{a}\cdot\dot{\Theta}(\partial_i,\partial_j)(0)\Bigr)\right.\right.
\\&\left.\left.\hspace{40pt}\cdot\exp\Bigl(\frac{1}{4}{\rm d}a{\rm
d}\overline{a}\cdot\dot{\Theta}(e_i,Je_i)(0)
-\frac{1}{2}F^{\mathcal{E}\!/\!S}(e_i,e_j)(0)e^i\wedge
e^j\Bigr)\right]\right)^{{\rm d}a{\rm d}\overline{a}}
\\
&=(2\pi \sqrt{-1})^{-n}\int_M\left(\widehat{\rm
A}\Bigl(\frac{1}{2}R_{klij}e^k\wedge e^l -\sqrt{-1}{\rm d}a{\rm
d}\overline{a}\cdot\dot{\Theta}(\partial_i,\partial_j)(0)\Bigr)\right.
\\
&\left.\hspace{20pt}\cdot\exp\Bigl(\frac{1}{4}{\rm d}a{\rm
d}\overline{a}\cdot\dot{\Theta}(e_i,Je_i)(0)
-\frac{1}{2}\tr_{T^{(1,0)}M}[R^+](0)\Bigr)\right)^{{\rm d}a{\rm
d}\overline{a}}\!\cdot\tr\bigl[\exp(-R^{\cal E})\bigr].
\end{split}
\end{align}

The facts that $\dot{\Theta}(e_i,Je_i)=2\tr_{T^{(1,0)}M}[U^+]$,
$\frac{1}{2}R_{ijkl}e^k\wedge e^l=-(R\partial_i,\partial_j)(0)$, and
$\dot{\Theta}(\partial_i,\partial_j)=\langle
UJ\partial_i,\partial_j\rangle$ imply that
\begin{align}
\begin{split}
&\widehat{\rm A}\Bigl(\frac{1}{2}R_{klij}e^k\wedge e^l
-\sqrt{-1}{\rm d}a{\rm
d}\overline{a}\cdot\dot{\Theta}(\partial_i,\partial_j)(0)\Bigr)
\\
&\hspace{30pt}\cdot\exp\Bigl(\frac{1}{4}{\rm d}a{\rm
d}\overline{a}\cdot\dot{\Theta}(e_i,Je_i)(0)
-\frac{1}{2}\tr_{T^{(1,0)}M}[R^+](0)\Bigr)
\\
&={\td}(R^+-{\rm d}a{\rm d}\overline{a}\,U^+).
\end{split}
\end{align}

Moreover, we get
\begin{equation}\label{83}
\begin{split}
&\lim_{t\rightarrow0}\left(\tr_s[e^{-I_t}]\right)^{{\rm d}a{\rm
d}\overline{a}}
\\
&=-(2\pi \sqrt{-1})^{-n}\int_M\left.\frac{\partial}{\partial
b}\right|_{b=0}\td(R^+ +b\,U^+)\cdot\tr\bigl[\exp(-R^{\cal E})\bigr].
\end{split}
\end{equation}

Combining Proposition \ref{81}, Proposition \ref{75}, Proposition
\ref{74}, (\ref{82}) and  (\ref{83}) together, we deduce
\begin{prop}\label{84}
The coefficient $a_0$ in Proposition \ref{81} can be
calculated by
\begin{equation}
a_0=-(2\pi\sqrt{-1})^{-n}\int_M\left.\frac{\partial}{\partial
b}\right|_{b=0}\td(R^+ +b\,U^+)\cdot\tr\bigl[\exp(-R^{\cal E})\bigr].
\end{equation}
\end{prop}

\

\makeatletter
\@startsection{subsection}{0}{0pt}{2mm}{0.2mm}{\bfseries}
{A proof of (\ref{j1})}\makeatother

Since the space of K\"{a}hler metrics
on $TM$ is convex, we may assume that $\ell\in\R\rightarrow g_{\ell}^{TM}$
is a smooth family of K\"{a}hler metrics on $TM$ such that $g_{0}^{TM}=g^{TM}$, $g_{1}^{TM}=g'^{TM}$.

Observing that the following algebraic identity holds,
\begin{equation*}
\tr^{\Omega^{p,\ast}(M)}\left[\exp\Bigl(\bigl\langle{\cal
A}\omega_i,\overline{\omega}_j\bigr\rangle\omega^i\wedge
i_{\omega_j}\Bigr)\right]=\sigma_p(\exp{\cal A}),\ \text{for any}\
{\cal A}\in \End(T^{1,0}M),
\end{equation*}
we get from Proposition \ref{81} and Proposition \ref{84} that
\begin{equation}
\begin{split}\label{j3}
M_{0,\ell}=&-(2\pi \sqrt{-1})^{-n}\!\int_M\frac{\partial}{\partial
b}_{|_{b=0}}\Biggl(\td(R_\ell^++b\,U_\ell^+)\cdot\sigma_p\Bigl(\exp(R_\ell^++bU_\ell^+)\Bigr)\Biggr)
\\
&\hspace{70pt}\cdot\tr\Bigl[\exp(-R^E)\Bigr].
\end{split}
\end{equation}

From Theorem \ref{59} and (\ref{j3}), we get

\begin{equation}\label{86}
\frac{\partial}{\partial\ell}\log\tau_{{\rm
holo},p,\ell}(M,E)
=(2\pi \sqrt{-1})^{-n}\!\int_M\frac{\partial}{\partial
b}_{|_{b=0}}\!\Bigl(\td_p(R_{\ell}^++b\,U_{\ell}^+)\Bigr)
\cdot\tr\Bigl[\exp(-R^E)\Bigr]\,.
\end{equation}

By the results of \cite[Section (e)]{bismut1988ata1}, the form
\begin{equation}
\varpi=(2\pi\sqrt{-1})^{-n}\!\int_{0}^1\frac{\partial}{\partial
b}_{|_{b=0}}\!\Bigl(\td_p(R_{\ell}^++b\,U_{\ell}^+)\Bigr){\rm d}\ell\cdot\tr\Bigl[\exp(-R^E)\Bigr]
\end{equation}
defines an element in $P/P'$ which depends only on $g^{TM}$ and $g'^{TM}$. According to \cite[Theorem 1.27, 1.29, and Corollary 1.30]{bismut1988ata1}, the component of degree $(n,n)$ of $\varpi$
represents in $P/P'$ the corresponding component of
\begin{equation*}
\widetilde{\td}_p\bigl(T^{(1,0)}M,g^{TM}, g'^{TM}\bigr)\cdot\ch\bigl(E,\nabla^E\bigr)\,.
\end{equation*}
Combining with (\ref{86}), we deduce that
\begin{equation}
\log\frac{\tau'_{{\rm holo},p}}{\tau_{{\rm holo},p}}=\int_M\varpi=\int_M\widetilde{\td}_p\bigl(T^{(1,0)}M,g^{TM}, g'^{TM}\bigr)\cdot\ch\bigl(E,\nabla^E\bigr),
\end{equation}
which is equivalent to (\ref{j1}).

The proof of Theorem \ref{j2} is completed.

\vspace{3mm}\th{Acknowledgements}%\dawuhao\.
The authors are indebted to Professor Weiping Zhang for his guidance
and very helpful discussions. 
%---------------------------------------------------------------参考文献
%\bibliographystyle{zhongguokexue}
%\bibliography{zhongguokexue}
\footnotesize
\makeatletter
\renewcommand\@biblabel[1]{\parindent=6mm
\REF{${#1}$\ }}
\makeatother
\def\refname{\no{\normalsize \bf References}}

\hml  %结束文稿的输入----------------------------------------------------------------------------------------------------

%-------------------------以下是可能要用到的命令
\begin{center}
\centerline{\psfig{figure=zkxf33.eps}} \centerline{\footnotesize
Fig. 1.\quad }
\end{center}

{\parbox[c]{60mm}\centerline{\psfig{figure=zkxf33.eps}}
\centerline{\footnotesize Fig. 1.\quad }} {}
{\parbox[c]{60mm}
\begin{center}
\footnotesize Table 1\quad \\\vspace{1.5mm} \doublerulesep 0.4pt
\tabcolsep 19pt
\begin{tabular}{\textwidth}{rcccc}
\hline \hline 表的内容 \hline \hline
\end{tabular}
\end{center}}
}

\newpage
\begin{center}
\footnotesize Table 1\quad \\\vspace{1.5mm} \doublerulesep 0.4pt
\tabcolsep 19pt
\begin{tabular*}{\textwidth}{rcccc}
\hline \hline 表的内容 \hline \hline
\end{tabular*}
\end{center}
\zihao{5}

*******************************************做图*********************
\begin{center}
\centerline{\psfig{figure=zkxf33.eps}} \centerline{\footnotesize
Fig. 1.\quad }
\end{center}

\parbox[c]{60mm}{\centerline{\psfig{figure=zkxf33.eps}}
\centerline{\footnotesize Fig. 1.\quad }}
\parbox[c]{60mm}
{ }

%双页码

%黑斜
\bf {\boldmath

%文字缩进
 \begin{enumerate}
 \item[(1)]
\end{enumerate}

%数字转罗文\romannumeral
%\公式居底如下
\begin{eqnarray}
 \dot {x} = Ax + Bu,  \quad x(0) = x_0,\nonumber\\
 \dot {x} = Ax + Bu,  \quad x(0) = x_0,
\end{eqnarray}

\def\no{\nonumber}\公式居底
\begin{eqnarray}
\end{eqnarray}

\begin{equation}公式齐缝
\end{equation}

\normalsize%每行居中

\underline%下划线